\documentclass[12pt]{article}
\usepackage{amssymb}
\usepackage{latexsym}
\usepackage{graphicx}
\usepackage{xspace}
\usepackage{placeins}
\usepackage{caption2}

\title{Computation in word-hyperbolic groups}
\author{David B. A. Epstein and Derek F. Holt}
\date{}

\newcommand{\Z}{\mathbb Z} 
\newcommand{\R}{\mathbb R} 
\newcommand{\N}{\mathbb N}
\newcommand{\Tau}{{\mathrm T}}
\newenvironment{proof}{\par\noindent\normalsize {\sc Proof}:}{{\hfill $\Box$}}
\newtheorem{theorem}{Theorem}[section]
\newtheorem{proposition}[theorem]{Proposition}
\newtheorem{lemma}[theorem]{Lemma}

\newcommand{\myvert}{\mbox{$\,\vert\,$}}
\newcommand{\mycolon}{\mbox{$\,:\,$}}
\newcommand{\intp}{meeting point\xspace}
\newcommand{\intv}{meeting vertex\xspace}
\newcommand{\intvs}{meeting vertices\xspace}
\newcommand\uast{^{\textstyle \ast}} 

\newcommand{\myfig}[2]{%
\begin{figure}[htbp]
{\centering \includegraphics{#1.ps}}
\caption{#2}
\label{#1}
\end{figure}
}
\newenvironment{mylist}{\begin{list}{}{
\setlength{\itemsep}{0mm}
\setlength{\parskip}{0mm}
\setlength{\topsep}{0mm}
\setlength{\parsep}{0mm}
\setlength{\itemsep}{0mm}
\setlength{\labelwidth}{7mm}
\setlength{\labelsep}{3mm}
\setlength{\itemindent}{0mm}
\setlength{\leftmargin}{10mm}
\setlength{\listparindent}{6mm}
}}{\end{list}}
\begin{document}
\maketitle
\section{Introduction} \label{introduction}
The purpose of this paper is to describe two algorithms for
computing with word-hyperbolic groups. Both of them have been implemented
in the second author's package {\sf KBMAG} \cite{KBMAG}.

The first is a method of
verifying that a group defined by a given finite presentation is
word-hyperbolic, using a criterion proved by Papasoglu in \cite{PAP},
which states that all geodesic triangles in the Cayley graph of a group
are thin if and only if all geodesic bigons are thin.
This is very similar to an algorithm described in \cite{WAK},
but it contains a simplification which appears to improve the performance
substantially. It also improves a less developed approach to the problem
described in \cite{HOLT}.

The second algorithm provides a method of estimating the constant of
hyperbolicity of the group with respect to the given generating set.
We do not know of any previously proposed general method for solving
this problem that has any prospect of being practical. Our current
implementation is experimental, and is very heavy on its use of
memory for all but the most straightforward examples, but it does at
least succeed on examples like the Von-Dyck triangle groups and the
two-dimensional surface groups.

Both of them follow a general philosophy of group-theoretical
algorithms that construct finite state automata. This approach
was originally proposed in the algorithms described in
Chapter 5 of \cite{ECHLPT} for computing automatic structures,
and employed in their implementations for short-lex structures
described in \cite{HOLT}. The basic idea is first to find
a method of constructing likely candidates for the required automata,
which we shall call the {\it working} automata,
The second step is to construct other (usually larger and more complicated)
{\it test} automata of which the sole purpose is to verify the correctness
of the working automata. In the case when this verification fails,
it should be possible to use words in the language of the test automata
to construct improved versions of the working automata.
One practical difficulty with this approach is that experience shows
that incorrect working automata and the resulting test automata
are much larger than the correct ones, so it can be extremely important
to find good candidates on the first pass.

The two algorithms dealt with in this paper are described in
Sections~\ref{verifying} and \ref{finding}. Some details
of their performance on some examples are presented in
Section~\ref{examples}. Unfortunately, we need to recall quite a lot of
notation from related earlier works, and we do this in Section~\ref{notation}.

There are a number of computational problems in which it is useful to
know hyperbolic constants which are either the same as those in this
paper, or are closely related to them. Some examples of such problems
follow.

In another paper we will describe a linear algorithm to
put a word in a word-hyperbolic group into short-lex normal form. The
linear estimate is due to Mike Shapiro (unpublished). Our method,
which is a bit different, may be usable in practice, though the ideas 
have not yet been implemented.
The standard algorithm for converting a word into normal
form in an automatic group is quadratic, as shown in \cite{ECHLPT}.

We also plan to show how to construct the automata which accept
1)~all bi-infinite geodesics, 2)~all pairs of asymptotic geodesics and
3)~all pairs of bi-infinite geodesics
which are within a finite Hausdorff distance of each other. Such
algorithms are necessary if one is to have any hope of a constructive
description of the limit space of a word-hyperbolic group, starting
with generators and relations of the group. Some of these automata
are also needed in Epstein's $n\log(n)$ solution of the conjugacy problem
(not yet published).

\section{Notation} \label{notation}
Throughout the paper, $G$ will denote a group with a given finite
generating set $X$. The identity element of $G$ will be denoted by $1_G$.
Let $A = X \cup X^{-1}$, and let $A\uast$ be the set of all words in $A$.
For $u,v \in A\uast$, we denote the image of $u$ in $G$ by $\overline {u}$,
and $u =_G v$ will mean the same as $\overline {u} = \overline {v}$.
For a word $u \in A\uast$, $l(u)$ will denote the length of $u$ and
$u(i)$ will denote the prefix of $u$ of length $i$, with $u(i) = u$
for $i \geq l(u)$.

Let $\Gamma = \Gamma_X(G)$ be the Cayley graph of $G$ with respect to $X$.
We make $\Gamma$ into a metric space in the standard manner, by letting
all edges have unit length, and defining the distance $\partial(x,y)$ between
any two points of $\Gamma$ to be the minimum length of paths connecting them.
(The points of $\Gamma$ include both the vertices, and points on the edges of
$\Gamma$.)
This makes $\Gamma$ into a {\it geodesic} space, which means that for
any $x,y \in \Gamma$ there exist geodesics (i.e. shortest paths) between
$x$ and $y$.
For $g \in G$, $l(g)$ will denote the length of a geodesic path
from the base vertex $1_G$ of $\Gamma$ to $g$.

A {\it geodesic triangle} in $\Gamma$ consists of three not
necessarily distinct points $a,b,c$ together with three directed geodesic
paths $u,v,w$ joining $bc$, $ca$ and $ab$, respectively.
The vertices $a,b,c$ of the triangle are not
necessarily vertices of $\Gamma$; they might lie in the interior of
an edge of $\Gamma$.

There are several equivalent definitions of word-hyperbolicity. The most
convenient for us is the following. Let $\Delta$ be a geodesic triangle
in $\Gamma$ with vertices $a,b,c$ and sides $u$, $v$, $w$ as above.
(Hence $l(u) = \partial(b,c)$, etc.) Let $\rho(a) = (l(v)+l(w)-l(u))/2$
and define $\rho(b), \rho(c)$ correspondingly. Then
$\rho(b)+\rho(c)=l(u)$, so any point $d$ on $u$ satisfies either
$\partial(d,b) \leq \rho(b)$ or $\partial(d,c) \leq \rho(c)$, and similarly
for $v$ and $w$. The points $d,e,f$ on $u,v,w$ with
$\partial(d,b) = \rho(b)$ and $\partial(d,c) = \rho(c)$, etc.,
are known as the {\it {\intp}s} of the triangle.

In a constant curvature geometry (the euclidean plane, the hyperbolic
plane or the sphere), the {\intp}s of a triangle are the points where
the inscribed circle meets the edges. In more general spaces, such as
Cayley graphs, the term {\it inscribed circle} has no meaning, but the
{\intp}s can still be defined.
\myfig{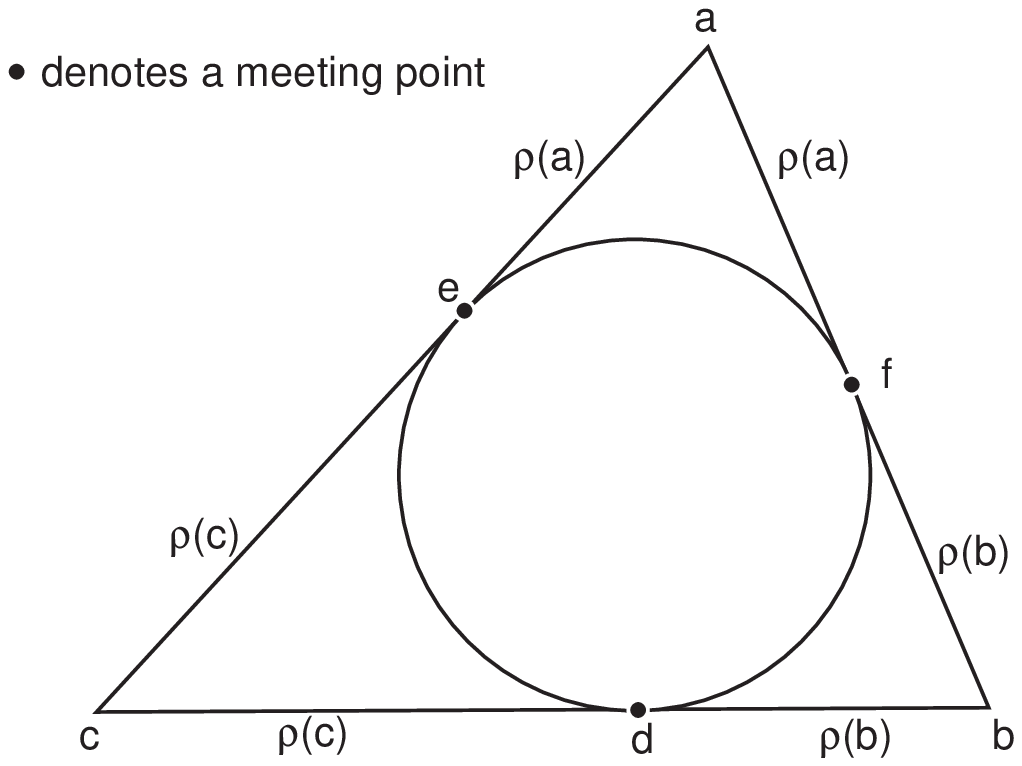}{This picture shows the {\intp}s as the intersection
of the inscribed circle with the edges of the triangle in the case of
constant curvature geometry.}

Suppose $a$, $b$ and $c$ are vertices in the Cayley graph.
Then the {\intp}s are also vertices if and only if
the perimeter $l(u) + l(v) + l(w)$ is even.

Let $\delta \in \R_+$. Then we say that $\Delta$ is {\it $\delta$-thin} if,
for any $r \in \R$ with $0 \leq r \leq \rho(x)$, the points $p$ and $q$ on
$v$ and $w$ with $\partial(p,a)=\partial(q,a)=r$ satisfy
$\partial(p,q) \leq \delta$, and similarly for the points within distance
$\rho(b)$ of $b$ and $\rho(c)$ of $c$. We call such points $p$ and $q$
{\it $\Delta$-companions}, or, if $\Delta$ is understood, just {\it companions}.
Note that the definition makes sense even when the triangle
$\Delta$ is not geodesic---we measure distances along the edges of the
triangle.
Normally companions are distinct, but there can be many situations
where they coincide---for example two geodesics sides of a triangle could have
an intersection consisting of a disjoint union of three intervals.
Mostly points on the triangle have exactly one companion, but the
{\intp}s normally have two companions---once again, in degenerate
situations two or all three of the {\intp}s may coincide.

The group $G$ is called word-hyperbolic if there exists a $\delta$ such
that all geodesic triangles in $\Gamma$ are $\delta$-thin. (It turns out
that this definition is independent of the generating set $X$ of $G$,
although the minimal value of $\delta$ does depend on $X$.)
The multi-author article \cite{ALO} is a good reference for the basic
properties of word-hyperbolic groups.

We also need to recall some terminology concerning finite state automata.
This has been chosen to comply with that used in \cite{HOLT} as far
as possible. The reader should consult \cite{ECHLPT} for more
details on the definitions and basic results relevant to the
use of finite state automata in combinatorial group theory.

Let $W$ be a finite state automaton with input alphabet $A$. We denote
the set of states of $W$ by $\mathcal{S}(W)$ and the set of initial and
accepting states by $\mathcal{I}(W)$ and $\mathcal{A}(W)$ respectively.
In a non-deterministic automaton there may be more than one transition
with a given source and label, and some transitions, known as
$\varepsilon$-transitions, may have no label.
In a deterministic automaton, there are no $\varepsilon$-transitions,
at most one transition with given source and label,
and $W$ has only one initial state.
(This type of automaton is named {\it partially deterministic}
in \cite{ECHLPT}.)
In this case, we denote the unique initial state by
$\sigma _0(X)$ and, for each $x \in A$ and
$\sigma \in \mathcal{S}(X)$, we denote the target of a
transition from $\sigma$ with label $x$ by $\sigma ^x$ if it exists.
We can then define $\sigma ^u$, for
$\sigma \in \mathcal{S}(X)$ and $u \in A\uast$, in the obvious way
whenever all of the required transitions exist.
The automata in this paper can be assumed to be deterministic unless
otherwise stated.

The automata that we consider may be one-variable or two-variable.
In the latter case, two words $u$ and $v$ in $A\uast$ are read simultaneously,
and at the same speed. This
creates a technical problem if $u$ and $v$ do not have the same length. To
get round this, we introduce an extra symbol $\$$, which maps onto the
identity element of $G$, and let $A^\dagger = A \cup \{\$\}$. Then if
$(u,v)$ is an ordered pair of words in $A$, we adjoin a sequence of
$\$$'s to the end of the shorter of $u$ and $v$ if necessary, to make them
both have the same length. The resulting pair will be denoted by
$(u,v)^\dagger$, and can be regarded as
an element of $(A^\dagger \times A^\dagger)\uast$. Such
a pair has the property that the symbol $\$$ occurs in at most one of $u$
and $v$, and only at the end of that word, and it is known as a {\it padded}
pair. We shall assume from now on,
without further comment, that all of the two-variable automata that arise
have input language $A^\dagger \times A^\dagger$ and accept only padded pairs.

Note that if $M$ is a two-variable automaton, then we can form a
non-deterministic automaton with language equal to
$\{u \myvert \exists v \mycolon (u,v)^\dagger \in L(M)\}$ simply by replacing
the label $(x,y)$ of each transition in $M$ by $x$.
(This results in an $\varepsilon$-transition in the case $x = \$$.)
We call this technique {\it quantifying over the second variable} of $M$.

Following \cite{HOLT}, we call a two-variable automaton $M$ a
{\it word-difference automaton} for the group $G$, if there is a function
$\alpha:\mathcal{S}(M) \rightarrow G$ such that
\begin{mylist}
\item[(i)] $\alpha(\sigma_0(M)) = 1_G$, and
\item[(ii)] for all $x,y \in A^\dagger$
and $\tau \in \mathcal{S}(M)$ such that $\tau ^{(x,y)}$ is defined, we have
$\alpha(\tau ^{(x,y)}) = x^{-1}\alpha(\tau)y$.
\end{mylist}
We shall assume that all states $\tau$ in a word-difference automaton $M$
are accessible; that is, there exist words $u,v$ in $A\uast$ such that
$\sigma _0(M)^{(u,v)^\dagger} = \tau$. It follows from
properties (i) and (ii) that $\alpha(\tau) = \overline {u^{-1}v}$,
and so the map
$\alpha$ is determined by the transitions of $M$.

Conversely, given a subset $\mathcal{D}$ of $G$ containing $1_G$,
we can construct a word-difference automaton $D$ with $\mathcal{S}(D) =
\mathcal{D}$, $\sigma_0(D)=1_G$ and, for $d,e \in \mathcal{D}$
and $x,y \in A^\dagger$ a transition $d \rightarrow e$ with label
$(x,y)$ whenever $x^{-1}dy =_G e$. The map $\alpha$ is the identity map.
We call this the word-difference automaton associated with $\mathcal{D}$.
(We have chosen not to specify the set $\mathcal{A}(D)$ of
accepting states of $D$, because this may depend on the context.)
Having constructed the automaton, we throw away the elements of
$\mathcal{D}$ which are not accessible from the initial state
$\sigma_0(D)$.

If $u,v \in A\uast$, then we call the set
$\mathcal{D} = \{\overline {u(i)}^{-1}\overline {v(i)} \myvert
i \in \Z, i \geq 0 \}$ the set of word-differences arising from
$(u,v)$. Then $(u,v)$ is in the language of the associated
word-difference machine $D$ provided that $\overline {u}^{-1}\overline {v}
\in \mathcal{A}(D)$.

The group $G$ is said to be automatic (with respect to
$X$), if it has an automatic structure.
This consists of a collection of finite state automata.
The first of these, denoted by $W$, is called the {\it word-acceptor}.
It has input alphabet $A$, and accepts at least one word in $A$ mapping
onto each $g \in G$.
The remaining automata $M_x$, are called the {\it multipliers}.
There is one of these for each generator $x \in A$, and also one for $x = 1_G$.
These are two-variable, and accept ($w_1,w_2)^\dagger$ for $w_1, w_2 \in A\uast$,
if and only if $w_1, w_2 \in L(W)$ and $w_1x =_G w_2$.
See \cite{ECHLPT} for an exposition of the basic properties of automatic
groups. It is proved in Theorem 2.3.4 of that book that there is a
natural construction of the multipliers $M_x$ of an automatic structure
as word-difference machines.

Now we fix a total order on the alphabet $A$.
The automatic structure is called {\it short-lex} if the language
$L(W)$ of the word-acceptor consists of the short-lex least
representatives of each element $g \in G$; that is the lexicographically
least among the shortest words in $A\uast$ that map onto $g$.
The existence of such a structure for a given group $G$ depends
in general on the generating set $X$ of $G$, but word-hyperbolic groups
are known to be short-lex automatic for any choice of generators. (This
is Theorem 3.4.5 of \cite{ECHLPT}.)

The group $G$ is called {\it strongly geodesically} automatic with respect
to $X$ if there is an automatic structure in which $L(W)$ is the set of all
geodesic words from $1_G$ to $g$ for $g \in G$. It is proved in Corollary
2.3 of \cite{PAP} that this is the case if and only if $G$ is
word-hyperbolic with respect to $X$ (from which it follows
that this property is independent of $X$). This result will be the basis
of our test for verification of word-hyperbolicity in Section~\ref{verifying};
the procedure that we describe will verify strong geodesic automaticity.

We shall assume throughout the paper that our group $G = \langle X
\rangle$ is short-lex automatic with respect to $X$, and that we have already
computed the corresponding short-lex automatic structure
$\{W, M_x \myvert x \in A \cup \{1_G\} \}$.
We assume also that the set $\mathcal{D}_M$ of all word-differences
that arise in the multipliers $M_x$ together with the associated
word-difference machine $D_M$ in which $1_G$ is the unique accepting state
has been computed.
These automata can be used to reduce (in quadratic time) words $u \in A\uast$
to their short-lex equivalent word in $G$ and so, in particular, we
can solve the word problem efficiently in $G$.
The above computations can all be carried out using the
{\sc KBMAG} package \cite{KBMAG}.

\section{Verifying hyperbolicity} \label{verifying}

Papasoglu (\cite{PAP}) has shown that a necessary and sufficient condition
for a group to be word-hyperbolic is as follows.
Let $\Gamma$ be the Cayley graph with respect
to some set of generators. The condition is that there is a number $c_P$, such
that, for any two geodesic paths $u,v:[0,\ell]\rightarrow \Gamma$ parametrised
by arclength, if $u(0)=v(0)$ and $u(\ell) = v(\ell)$, then, for all $t$
satisfying $0\le t \le \ell$, $d_\Gamma(u(t),v(t)) \le c_P$.
The least possible value of $c_P$ is called
{\it Papasoglu's constant}. Such a configuration of
$u$ and $v$ is called a {\it geodesic bigon}.
In order to know that all geodesic bigons have
uniformly bounded width, there is no
loss of generality in restricting
to the case where $u(0)=v(0) $ is a vertex of the Cayley graph. We can
also restrict to the case where
$u(\ell)=v(\ell)$ is either a vertex of the Cayley graph
or the midpoint of an edge.
It is unknown whether the uniform thinness of such geodesic bigons
follows from the uniform thinness of the more special
geodesic bigons with both ends vertices.

Our algorithm not only verifies word-hyperbolicity for a given group,
but also gives a precise computation of the smallest possible value of
Papasoglu's constant. In fact, it gives even more precise information, namely
the set of word-differences $\overline{u(i)}^{-1}\overline{v(i)}$,
where $u$ and $v$ vary over all geodesic bigons with $u(0) = v(0)$ a vertex
and $i$ varies over all positive integers. In all the examples we have looked
at, we have observed that the number
of such group elements is very much smaller than the number of group elements
of length at most $c_P$. This is important in practical computations.

The algorithm for verifying word-hyperbolicity proceeds by
constructing sequences $WD_n$, $ GE_n$,
$GW_n$, $T_n$, for $(n > 0)$, of finite state automata.
(These letters stand respectively for `word-difference', `geodesic-equality',
`geodesic word acceptor' and `test'.)
In general, for $n>0$, $\mathcal{WD}_n$ will be a set of elements of $G$
containing $\{1_G\}$,
and $WD_n$ will be the associated word-difference machine in which
$1_G$ is the only accepting state.

We define $\mathcal{WD}_1$ to be the
set $\mathcal{D}_M$ defined at the end of Section~\ref{notation},
and let $W$ be the short-lex word-acceptor.
Then, when $n \geq 1$, we define $GE_n$, $GW_n$ and $T_n$ as follows. We
define the language
$$L(GE_n) =
\{(u,v) \in A\uast \times A\uast \myvert (u,v) \in L(WD_n), v \in L(W), l(u)=l(v)\}.$$
Recall that all elements of $L(W)$ are geodesics and that
the only accept state of $WD_n$ is $1_G$.
It follows that in the previous definition, $\overline u = \overline v$
and that $u$ and $v$ are both geodesics.
Now we define the language
$$L(GW_n) = \{u \in A\uast \myvert
\exists v \in A\uast \mycolon (u,v) \in L(GE_n) \}.$$
Again, $u$ must be a geodesic.

Then we define the language
\begin{eqnarray*}
\lefteqn{L(T_n) =
\{w \in A\uast \setminus L(GW_n) \myvert}\\
&&\exists u\mycolon
(w,u) \in L(WD_n), u \in L(GW_n), l(u)=l(w)\}.
\end{eqnarray*}
Again, $u$ and $w$ are both geodesics in the previous definition.

\myfig{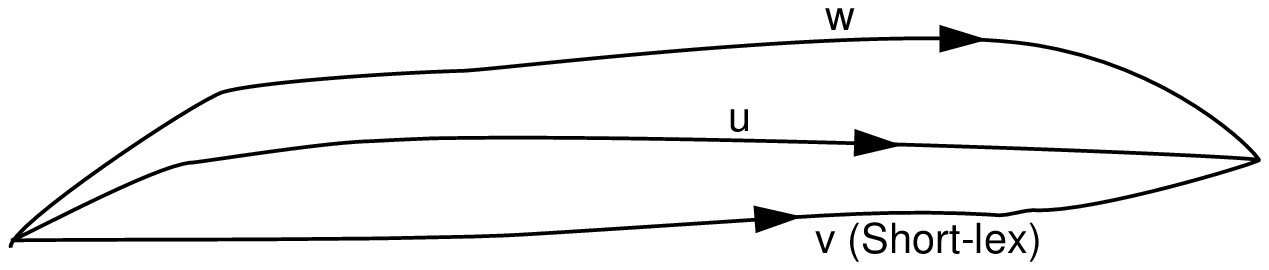}{This illustrates the geodesic paths $u$, $v$ and $w$
described in Section~\ref{verifying}.}

If $L(T_n)$ is empty for some $n$, then the procedure halts.
Otherwise, we find a geodesic word $w \in L(T_n)$, reduce it to its
short-lex least representative $v$, and define $\mathcal{WD}_{n+1}$
to be the union of $\mathcal{WD}_n$ and the set of word-differences
arising from $(w,v)$. Then we can define the automaton
$WD_{n+1}$ and construct the other
automata for the next value of $n$.

\begin{theorem} The above procedure halts if and only if $G$ is
strongly geodesically automatic with respect to $X$.
\end{theorem}
\begin{proof}
First note that, if $L(T_n)$ is non-empty for some $n$,
and contains the word $w$ reducing to $v \in L(W)$,
then the word-differences arising from $(w,v)$ cannot all lie in
$\mathcal{WD}_n$. For otherwise we would have $(w,v) \in L(WD_n)$, and
hence $(w,v) \in L(GE_n)$ and $w \in L(GW_n)$. But $w$ has been chosen
so that it is not in $L(GW_n)$. Hence
$\mathcal{WD}_{n+1}$ strictly contains $\mathcal{WD}_n$. Thus, if
the procedure does not halt, then the number of word-differences arising
from pairs of geodesics $(u,v)$ with $u =_G v$ cannot be finite,
and so by Theorem 2.3.5 of \cite{ECHLPT}
$G$ cannot be strongly geodesically automatic.

Conversely, suppose that the procedure does halt, and that $L(T_n)$
is empty for some $n \geq 1$. We claim that $L(GW_n)$ is equal to the
set of all geodesic words. If we can show this, then it would follow from the
definition that $L(GE_n)$ is equal to the set of pairs of
geodesic words $(u,v)$ such that $v \in L(W)$ and $u =_G v$.
But then, if $u_1,u_2 \in A\uast$ are geodesics with
$l(\overline{u_1}^{-1}\overline{u_2}) \leq 1$ and $u_1,u_2$
reduce to $v_1,v_2 \in L(W)$, then $(u_1,v_1), (u_2,v_2) \in
L(WD_n)$ whereas $(v_1,v_2) \in L(D_M)$.
This would show that the word-differences
arising from $(u,v)$ have bounded length, and so $G$ is strongly geodesically
automatic by Theorem 2.3.5 of \cite{ECHLPT}. Thus it suffices to establish
the claim that $L(GW_n)$ is equal to the set of all geodesic words.

Suppose that the claim is false, and let $u$ be a minimal length
geodesic word not contained in $L(GW_n)$. Since $L(GW_n)$ contains
the empty word,
$u=wx$ for some word $w$ and some $x \in A$.
Then $w$ is a geodesic word and $l(w) = l(u)-1$.
Suppose that $w$ reduces to $v \in L(W)$. Then $w \in L(GW_n)$
and $(w,v) \in L(GE_n)$. So $(w,v) \in L(WD_n)$.
Therefore $(u,vx) = (wx,vx) \in L(WD_n)$.

If $vx \in L(W)$, then by definition of $GE_n$ we get
$(u,vx) \in L(GE_n)$, and so $u \in L(GW_n)$, a contradiction.
On the other hand, if $vx \notin L(W)$ and $vx$ reduces to $v' \in L(W)$,
then $(v,v') \in L(M_x)$. From the definition of $WD_1$
we have $(vx,v') \in L(GE_1)$,
which implies $vx \in L(GW_r)$ for all $r \geq 1$.
But then $u \in L(T_n)$, contrary to assumption.
\end{proof}

The above procedure is very similar to that described in
\cite{WAK}. The principal difference is that our definition of the
test-machines $T_n$ is rather simpler. Furthermore,
in our implementation, we do not construct the non-deterministic
automaton resulting from the quantification over the second variable in
the definition of $L(T_n)$. Instead, we construct the two-variable
automaton with language
$$ \{(w,u) \myvert 
\ (w,u) \in L(WD_n),\ u \in L(GW_n), \ l(u)=l(w)\},$$
and check during the course of the construction whether
there are any words $w \notin L(GW_n)$ that arise.
If the construction completes and we find
no such $w$, then $L(T_n)$ is empty.
Otherwise we abort after finding some fixed number of words (such as 500)
$w \notin L(GW_n)$ and use all of these words with their short-lex
reductions to generate new word-differences.

So the procedure stops if and only if $G$ is word-hyperbolic. If it stops
at the $n$-th stage,
then $\mathcal{WD}_n$ is a finite set of word-differences which gives
the best possible value of Papasoglu's constant for the particular choice
of generators.
More precisely, the procedures described above find
the set of all word-differences for all pairs of geodesics $(w,u)$
which start at the same vertex of the Cayley graph and
end at vertices at distance at most one apart,
and where the word $u$ is short-lex minimal.
But by using a standard composite operation on two-variable automata
as defined in~\cite{HOLT} for example,
we can easily compute from this the word-differences for general
geodesic bigons in which at least one of the vertices is a vertex of the
Cayley graph.
Moreover, we have constructed an automaton whose language is the set
of all geodesics in the Cayley graph that start and end at vertices
of the Cayley graph.

\section{Finding the constant of hyperbolicity} \label{finding}
Throughout this section, we assume that $G = \langle X \rangle$ is a
word-hyperbolic group and that $\delta>0$ is a constant such that
all geodesic triangles in $\Gamma_X(G)$
are $\delta$-thin. The aim is to devise a practical
algorithm to find such a $\delta$, which should of course be as small as
possible. As before, we assume that we have already calculated the short-lex
automatic structure for $G$ with respect to $X$.

\subsection{The reverse of a finite state automaton}
\label {reverse}
Our procedure makes use of reversed automata, and so we start with
a brief discussion of this topic.

If $w \in A\uast$ is a word, then we denote the reversed word by $w^R$.
Let $M$ be a finite state automaton with alphabet $A$. We want to
form the reversed automaton $M^R$ with language $\{w^R \myvert w \in L(M)\}$.
We can define a non-deterministic version $NM^R$ of $M^R$, simply by
reversing the arrows of all transitions, and interchanging the sets of initial
and accepting states. Then we can build a deterministic machine $M^R$
with the same language in the standard way, by
replacing the set $\Sigma = \mathcal{S}(NM^R)$ of states of $NM^R$ with its
power set $\mathcal{P}(\Sigma)$, and, for $\Tau,\Upsilon \in \mathcal{
P}(S)$, defining a transition in $M^R$ with a label $x$ from $\Tau$ to
$\Upsilon$, if $\Upsilon$ is the set of all states $\upsilon$ of $M$
from which there exists an arrow labelled $x$ in $M$ with target some
$\tau\in \Tau$.
If we think of $x$ as a partial map $p_x$ from $\mathcal S (M)$ to itself,
then the existence of an arrow in $M^R$ from $\Tau$ to $\Upsilon$ is
equivalent to saying that
$\Upsilon$ is the full inverse image of $\Tau$ in $\mathcal S
(M)$ under $p_x$.

The initial state of $M^R$ is the set of all accepting states of $M$ and
a state of $M^R$ is accepting whenever it contains an initial state of $M$.

In practice, we do not need to use the complete power set
$\mathcal{P}(\Sigma)$. We start with the set of accepting states of $M$ as
initial state of $M^R$, and then construct the accessible states and
transitions of $M^R$ as the orbit of the initial state under the action of
$A\uast$.

Let $\sigma_0(M)$ and $\sigma_0(M^R)$ be the initial states of $M$ and
$M^R$ respectively.
Suppose that $v,w \in A\uast$, $\sigma_0(M)^v = \tau$, and $\sigma_0(M^R)^w =
\Tau \subseteq \mathcal{P}(\Sigma)$. Then, from the construction above, we
see that $\tau \in \Tau$ if and only if $vw^R \in L(M)$. We shall need
this property below, and when we compute the reverse of an automaton
we need to remember the subsets of $\Sigma$ that define the states of
$M^R$. (This means that we cannot minimise $M^R$, but in practice this
does not appear to be a problem because, at least for word-acceptors of
automatic group, $M^R$ does not seem to have many more states than $M$.)

\subsection{Reduction to short-lex geodesic triangles}
The reader should now recall the notation for $\delta$-thin hyperbolic
triangles with vertices $a,b,c$ and edges $u,v,w$ in the Cayley graph
$\Gamma = \Gamma_X(G)$ defined in Section~\ref{notation}.
\myfig{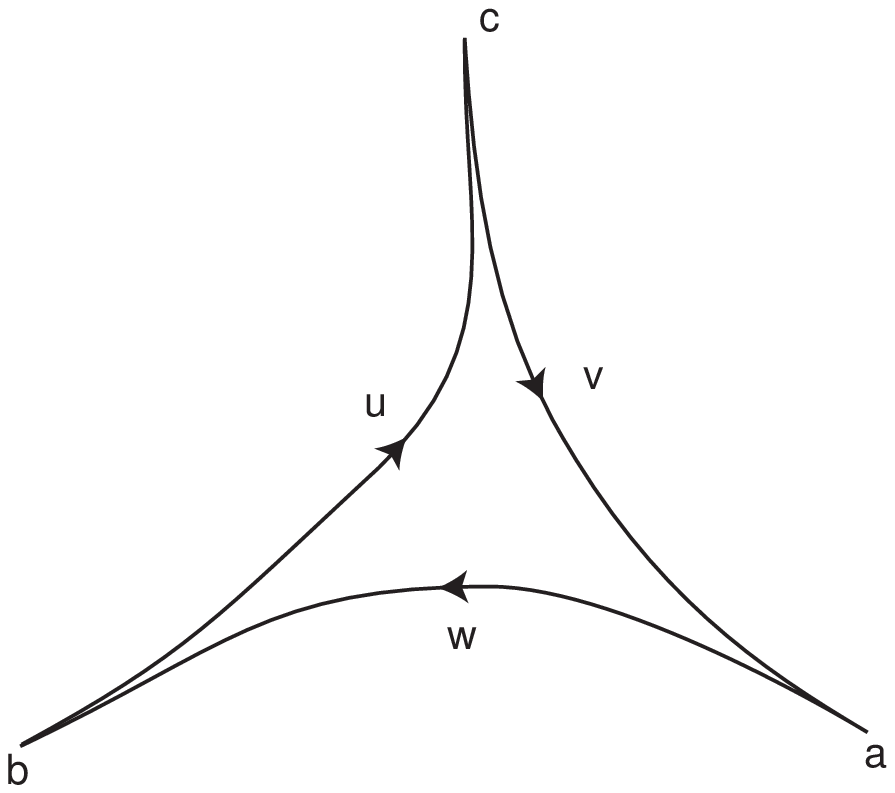}{A triangle with short-lex sides, annotated as in the
text.}
We shall call a geodesic
triangle in $\Gamma$ {\it short-lex geodesic}, if its vertices are vertices
of $\Gamma$ and if the words $A\uast$ corresponding to the
edges of the triangle (which we shall also denote by $u,v,w$)
all lie in $L(W)$; that is, they are all short-lex minimal words.

It is important to work with short-lex triangles, because in general
there are far more geodesic triangles and consideration of all of
these is likely to make an already difficult computational problem
impossible.

Our algorithms are designed to compute the minimal $\delta \in \N$
such that all short-lex geodesic triangles are $\delta$-thin.
In fact they do considerably more than this, because they compute
the set of word-differences that arise from the word-pairs $(u,v)$,
where $u$ and $v$ are the words from the two sides of such a triangle
that go from a vertex of the triangle as far as the {\intp}s
on the two sides.

From the following proposition, we derive a bound on the
value of the thinness constant of hyperbolicity for arbitrary geodesic
triangles. As before, let $\mathcal{D}_M$ be the set of word-differences
associated with the multiplier automata $M_x$ in the short-lex automatic
structure of $G$, and let $\gamma$ be the length of the longest element in
$\mathcal{D}_M$. Let $\gamma'$ be the length of the longest element in the
final stable set $\mathcal{WD}_n$ of word-differences defined in
Section~\ref{verifying}. Although the bound in the next
proposition will in general be too large, the fact that $\gamma$ and $\gamma'$
are usually smaller than $\delta$ means that it is likely only to
be wrong by a small constant factor.

\begin{proposition}\label{short-lex thin implies thin}
Suppose that all short-lex geodesic triangles in $\Gamma$ are $\delta$-thin.
Then all geodesic triangles in $\Gamma$ are
$(\delta+2(\gamma+\gamma')+3)$-thin.
\end{proposition}
\begin{proof}
Let $\Delta$ be any geodesic triangle in $\Gamma$. We fix a direction
around the triangle which we call {\it clockwise}. We define a new
triangle $\Delta'$, which is equal to $\Delta$ except that its vertices
are (if necessary) moved clockwise around $\Delta$ to the nearest vertex
of $\Gamma$. Thus each vertex is moved through a distance less than one 1,
and it follows easily that the {\intp}s of the triangle move through
a distance less than 2. It is not difficult to see that if $p$
and $q$ are $\Delta$-companions (see Section~\ref{notation}
for definition), then there exist $\Delta'$-companions $p'$ and $q'$ with
$\partial(p,p') + \partial(q,q') < 2$. (A careful argument is needed
if $p$ and $q$ are near {\intp}s.) It follows that, if $\Delta'$ is
$\delta'$-thin for some $\delta'$, then $\Delta$ is ($\delta'+2$)-thin.

Let $a,b,c$ be the vertices of $\Delta'$ in clockwise order.
Let $a'$ be the vertex on the union of the three sides of $\Delta'$,
adjacent to $a$ in an anticlockwise direction, and similarly for $b'$ and $c'$.
(So the vertices of
$\Delta$ lie on the edges $(a',a), (b',b)$ and $(c',c)$.) We will assume
that the six vertices $\{a,a',b,b',c,c\}$
are distinct, and leave to the reader the minor modifications
needed for the cases where two or more of them coincide.

Let $u$, $v$ and $w$ be the paths from $b$ to $c$, $c$ to $a$ and $a$ to $b$,
respectively, defining short-lex reduced words $u$, $v$ and $w$ in $L(W)$,
and let $\Delta''$ be the resulting short-lex geodesic triangle with
vertices $a$, $b$ and $c$.
The triangle $\Delta'$ is not necessarily geodesic, but
the paths $u',v',w'$ on $\Delta'$ from $b$ to $c'$, $c$ to $a'$ and
$a$ to $b'$, respectively,
are geodesic, because they are subpaths of the sides of $\Delta$.

Let $u''$, $v''$ and $w''$ be the paths from $b$ to $c'$, $c$ to $a'$ and
$a$ to $b'$, respectively, defining short-lex reduced words in $L(W)$.
The situation is illustrated in Figure~\ref{sltriangle}.
\myfig{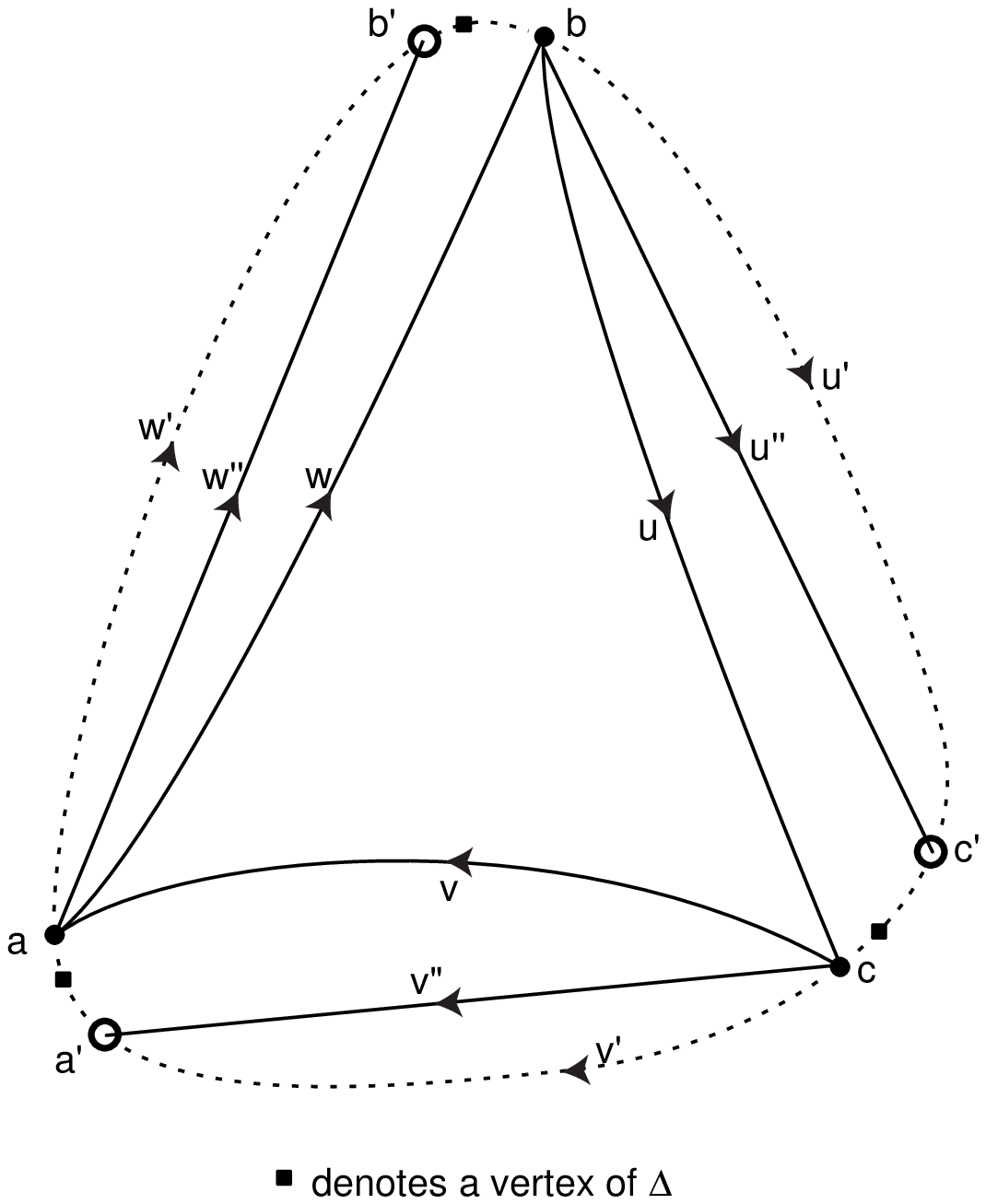}{This diagram illustrate the situation described in the
proof of Proposition~\ref{short-lex thin implies thin}.
}

Then $u$ and $u''$ are words in $L(W)$ that have the same starting point in
$\Gamma$ and end a distance one apart. So the word-differences arising
from $(u,u'')$ lie in $\mathcal{D}_M$. Since the paths $u'$ and $u''$
have the same starting and end points, $u'$ is a geodesic and $u'' \in L(W)$,
the word-differences arising from $(u',u'')$ lie in $\mathcal{WD}_n$.
It follows that the word-differences arising from $(u,u')$ have length
at most $\gamma + \gamma'$, and similarly for $(v,v')$ and $(w,w')$.

\myfig{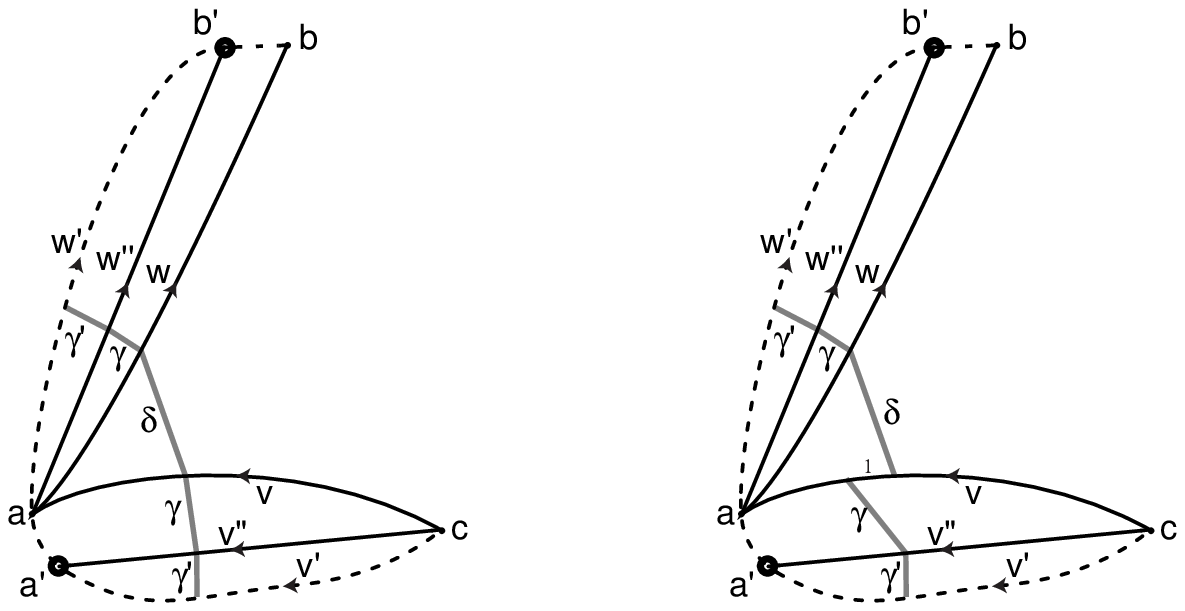}{These illustrates the origin of $\delta + 2(\gamma+\gamma')+1$
in the expression bounding $\partial(p',q')$.
In the diagram on the left, $l(ca) = l(ca') + 1$, and in the
diagram on the right, $l(ca) = l(ca')$.
The shapes of the curves are due to the fact that we are sometimes
looking at points which are equidistant from $a$ and sometimes at
points which are equidistant from $c$.
Interested readers are left to work out the details for themselves.
}

The triangles $\Delta'$ and $\Delta''$ have the same vertices $a$, $b$ and
$c$.
Each side of $\Delta'$ is at least as long as the corresponding side of
$\Delta''$, but no more than one unit longer---this can be deduced from the
fact that the sides of the original triangle $\Delta$ are geodesic.
It follows that the distance from a vertex
of $\Delta'$ to its two adjacent $\Delta'$-{\intp}s is within one unit of its
distance to its two adjacent $\Delta''$-{\intp}s.
We then see that, if all pairs $p,q$ of $\Delta''$-companions
satisfy $\partial(p,q) \leq \delta$, then all pairs
$p',q'$ of $\Delta'$-companions satisfy
$\partial(p',q') \leq \delta + 2(\gamma+\gamma')+1$.
The result now follows.

Part of the argument is illustrated in Figure~\ref{corner}
\end{proof}

\subsection{The automata $FRD$ and $FRD^3$ }
In this section, we describe a finite state automaton which we shall call
$FRD$, which stands for `forward, reverse, difference'. Roughly speaking,
it is a two-variable machine which reads the two sides emerging from a vertex
of a short-lex geodesic triangle as far as the {\intp}s on those two
sides. We also describe an associated automaton $FRD^3$ which consists of
three copies of $FRD$.
The three pairs of words read by $FRD^3$ will be accepted when they are the
three pairs of edges emerging from the three vertices of some short-lex geodesic
triangle, ending and meeting at the {\intp}s of the triangle.

\myfig{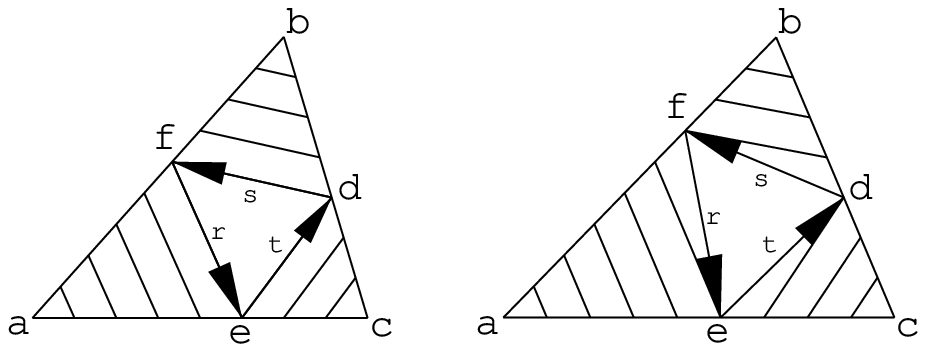}{This diagram shows the \intvs, denoted by $d$, $e$ and $f$.
Companions have been joined by a line.
The diagram on the left shows the situation where the perimeter of the
triangle is even and the diagram on the right where the perimeter is odd.}

When the perimeter $l(u)+l(v)+l(w)$ of a short-lex geodesic triangle is even,
the {\intp}s $d,e,f$ are vertices of $\Gamma$ which lie on $u=bc$, $v=ca$
and $w=ab$ respectively.
When the perimeter is odd, however, each {\intp} lies in the middle of
an edge of $\Gamma$, which is not convenient for us.
We therefore move them to a neighbouring vertex and re-define $d,e,f$ to
be the points on $u,v,w$ at distance $\rho(b)+1/2$, $\rho(c)+1/2$ and
$\rho(a)+1/2$ from $b,c$ and $a$, respectively, and call
$d,e,f$ the {\it \intvs}.
Each of the three vertices of the triangle has two incident sides, one
in the clockwise direction and the other in the anti-clockwise direction,
starting from the vertex. In these terms, if we start
at a vertex of the triangle and move away from this vertex along the two sides
of the triangle emerging from it, one edge at a time,
then we have to move one extra edge along the clockwise side than the
anti-clockwise side to reach the \intvs on the two sides.

Let $W$ be the word acceptor for the short-lex automatic structure of $G$.
We assume that the reverse $W^R$ of $W$
has been computed, as described in~\ref{reverse}.
For a short-lex geodesic triangle
with \intvs $d$, $e$ and $f$ defined as above, we
denote the elements of $G$ corresponding to paths from $d$ to $f$, $e$ to $d$,
and $f$ to $e$, by $r$, $s$ and $t$, respectively. This is illustrated in
Figure~\ref{intv}.

Let $\mathcal{D}_1$ be the
set of elements $r$ of $G$ that arise by considering all such triangles.
(By symmetry, this set includes all of the elements $s$ and $t$ as
well.)
Let $\mathcal{D}_2$ denote the set of all elements of $G$ of the
form
$$\{\overline {w(i)}^{-1}\overline {v^R(i)} \myvert i \in \Z, 0 \leq i \leq
\rho(a)\},$$
for all triangles under consideration; that is, the set of word-differences
arising from reading the two edges of the triangle from $a$ up to
the {\intp}s. (For triangles with even perimeter,
the elements $r,s,t$ lie in both $\mathcal{D}_1$ and $\mathcal{D}_2$.)
Let $\mathcal{D}_T = \mathcal{D}_1 \cup \mathcal{D}_2$. Then our
assumption that geodesic triangles are $\delta$-thin implies that
$\mathcal{D}_T$ is finite, and that its elements have length at
most $\delta + 1$.

The automaton $FRD$ is defined as follows. Its states are triples
$(\sigma,\Sigma,g)$ with $\sigma \in \mathcal{S}(W)$, $\Sigma \in
\mathcal{S}(W^R)$ and $g \in \mathcal{D}_T$. (This explains the name $FRD$.)
For $x \in A$ and $y \in A^\dagger$,
a transition is defined with label $(x,y)$ from
$(\sigma,\Sigma,g)$ to $(\sigma',\Sigma',g')$ if and only if
$\sigma^x = \sigma'$, $\Sigma^y = \Sigma'$ and $x^{-1}gy = _G g'$.
(We allow $y$ but not $x$ to be the padding symbol because, in the case
of triangles with odd perimeter, we need to read one extra generator from
the left hand edge than the right hand edge to get to the \intv.
We define $\Sigma^\$ = \Sigma$ for all $\Sigma$.)
The initial state is $(\sigma_0(W),\sigma_0(W^R),1_G)$.
The accepting states of $FRD$ are those with the third component
$g \in \mathcal{D}_1$.

The above description is not quite correct. In fact, each state has a fourth
component which is either 0 or 1, and is 1 only when a pair $(x,\$)$ has
been read. There are no transitions from such states.

We also need to consider an automaton $FRD^3$, which consists of
the product of
three independent copies $FRD_a$, $FRD_b$ and $FRD_c$ of $FRD$.
Its input consists of sextuples of
words $(u_a,v_a,u_b,v_b,u_c,v_c) \in (A\uast)^6$, where $(u_a,v_a)$,
$(u_b,v_b)$ and $(u_c,v_c)$ are input to the three copies
of $FRD$. A state of $FRD^3$ consists of a triple $(\tau_a,\tau_b,\tau_c)$,
where $\tau_a = (\sigma_a,\Sigma_a,g_a)$ is a state of $FRD_a$, etc.
The initial state of $FRD^3$ is the triple consisting of the initial
states of the three copies.

To specify the set $\mathcal{A}(FRD^3)$ of accepting states,
we recall that a state $\Sigma$ of $W^R$ is a
subset of the set $\mathcal{S}(W)$ of states of $W$. We have
$(\tau_a,\tau_b,\tau_c) \in \mathcal{A}(FRD^3)$ if and only if
$\sigma_a \in \Sigma_b$, $\sigma_b \in \Sigma_c$, $\sigma_c \in \Sigma_a$,
and $g_cg_bg_a = 1_G$.

The proof of the following lemma should now be clear.
\begin{lemma} The triple $(u_a,v_a,u_b,v_b,u_c,v_c)$ is accepted by $FRD^3$ if
and only if $u_bv_c^R$, $u_cv_a^R$ and $u_av_b^R$ all lie in $L(W)$ and
form the three sides of a short-lex geodesic triangle in $\Gamma$.
\end{lemma}

\subsection{Proving correctness of $FRD$}
If we can construct $FRD$, then we can compute the value of the hyperbolic
thinness constant $\delta$ for short-lex geodesic triangles
as the maximum of the lengths of the words in $\mathcal{D}_T$.
(Or more exactly, this number plus 1, because of our re-definition of
the {\intp}s of the triangles.)
Conversely, we have already computed $W$ and $W^R$, so the construction
of $FRD$ only requires us to find the set $\mathcal{D}_T$.
The idea is to construct a candidate for $FRD$, and then
to verify that it is correct.

In one of our short-lex geodesic triangles, two of the edge-words
$u,v$ can be chosen as arbitrary elements of $L(W)$, and the third is
then determined as the representative in $L(W)$ of
$(\overline {uv})^{-1}$. So we can proceed by choosing
a large number of random pairs of words in $u,v \in L(W)$
(for example we might choose 10000 such pairs of length up to 50),
computing $w$ as described, and then computing the set $\mathcal{D}_T$
of word-differences arising from all of these triangles. We do this
several times until the set $\mathcal{D}_T$ appears to have stopped
growing in size. (We are assuming that $G$ is word-hyperbolic, so we
know that the correct set $\mathcal{D}_T$ is finite.)
We then proceed to the verification stage, which is computationally
the most intensive.

The idea is to compute a two-variable finite state automaton $GP$ (geodesic
pairs) of which the
language is the subset of $A\uast \times A\uast$ defined by the expression
\begin{eqnarray*}
\lefteqn{\{(w_1,w_2)^\dagger \in A\uast \times A\uast \myvert}\\
&&\exists (u_a,v_a,u_b,v_b,u_c,v_c) \in L(FRD^3) \mycolon w_1 = u_av_b^R,
w_2 = v_au_c^R \}.
\end{eqnarray*}
Then $GP$ accepts the set of pairs of sides $(w,v^R)$ emerging from the
vertex $a$ in the triangles that are accepted by the current version of
$FRD^3$. (During the course of the program, $FRD$ changes as more
information is incorporated into it.)
Thus $FRD$ is correct if and only if $L(GP) = L(W) \times L(W^R)$.
Since checking for equality of the languages of minimised deterministic
automata is easy, we can perform this check provided that we can construct
$GP$.

Furthermore, if the check fails then our definition of
$GP$ ensures that $L(GP) \subset L(W) \times L(W^R)$. So we can find one
or more specific words $ (w_1,w_2) \in L(W) \times L(W^R) \setminus L(GP)$
and then compute the word-differences arising from the short-lex geodesic
triangle having $w_1$ and $w_2^R$ as two of its sides.
We can then adjoin these to $\mathcal{D}_T$ and return to the
construction of $FRD$.

The construction of $GP$ can be carried out in principle, but because of the
large number of quantified variables involved in the above expression, a naive
implementation would be hopelessly expensive in memory usage.

We now give a second and more detailed version of our implementation of
this construction, in a way that makes the computation easier to carry
out. It remains heavy in its memory usage, but it does at
least work for easy examples.
The basic objects collected during the course of the computation
are the word-differences $\mathcal{D}_T$
referred to above, which are used in constructing $FRD$,
and the triples $(r,s,t)$ of small triangles in the Cayley
graph whose vertices are the \intvs of some short-lex geodesic triangles.
The situation is shown in Figure~\ref{intv}.

The general idea is to map out a short-lex geodesic triangle by
advancing from one vertex, say $a$, using $FRD$.
Then, at a certain moment, the machine `does the splits', with a
non-deterministic jump corresponding to one of the small $(r,s,t)$-triangles.
To complete the triangle, the two legs have to be followed to the
other vertices $b$ and $c$, respectively, of the large triangle. Each of
the legs follows the reverse of $FRD$.

Here is a third, even more explicit, version of the construction,
which can be skipped by readers who are only interested in the
conceptual description of the program.
We assume that we have constructed a candidate
for $FRD$ explicitly. We do not construct $FRD^3$ explicitly, but we do make a
list of all triples $(r,s,t)$ such that $r,s,t \in \mathcal{S}(FRD)$
and $(r,s,t) \in \mathcal{A}(FRD^3)$.
Having done that, we can forget the structure of the states of $FRD$ as
triples, and simply manipulate them as integers.

We also need a version of an automaton $FRD^R$ that accepts the reverse
language of $FRD$. In this case it is convenient to work with a partially
non-deterministic version---that is, it is deterministic except that
there are many initial states, not just one.
The states are subsets of $\mathcal{S}(FRD)$
as described in \ref{reverse} and the transitions and
accepting states are also as described there.
But instead of having
a unique initial state, for each accepting state $\tau$ of $FRD$ we
make the singleton subset $\{\tau\}$ into an initial state of $FRD^R$.
Note also that if $(u,v)^\dagger \in L(FRD)$ with $l(u)=l(v)+1$, then
the reversed pair accepted by $FRD^R$ is $(u^R,\$v^R)$; that is, the
padding symbol comes at the beginning rather than the end.  In general,
let us denote the word-pair formed by inserting the padding symbol at
the beginning by $^\dagger(u,v)$.

We shall now describe a non-deterministic version $NGP$ of $GP$
(geodesic pairs).
Subsequent to its construction,
it can be determinised using the usual subset
construction, minimised, and its language compared with $L(W) \times L(W)^R$.

The precise description of $NGP$ is rather technical, so
we shall first attempt to explain its operation. A (padded) pair of words
$(w_1,w_2)^\dagger$ is to be accepted if and only it satisfies the expression
above; that is, $(w_1,w_2) = (u_av_b^R,v_au_c^R)$, where $(u_a,v_a)^\dagger
\in L(FRD)$, and there exist words $u_b,v_c \in A\uast$ with
$(u_b,v_b)^\dagger, (u_c,v_c)^\dagger \in L(FRD)$
and $(u_a,v_a,u_b,v_b,u_c,v_c) \in L(FRD^3)$.

Equivalently, (writing the reverses of $u_b$, $u_c$, $v_b$ and $v_c$
by capitalising $u$ and $v$), $(w_1,w_2)^\dagger$
is accepted if and only of $(w_1,w_2) = (u_aV_b,v_aU_c)$ where
$(u_a,v_a)^\dagger \in L(FRD)$, and there exist words $U_b$, $V_c \in A\uast$
with $^\dagger(U_b,V_b)$, $^\dagger(U_c,V_c) \in L(FRD^R)$
and $$(u_a,v_a,U_b^R,V_b^R,U_c^R,V_c^R) \in L(FRD^3).$$

\myfig{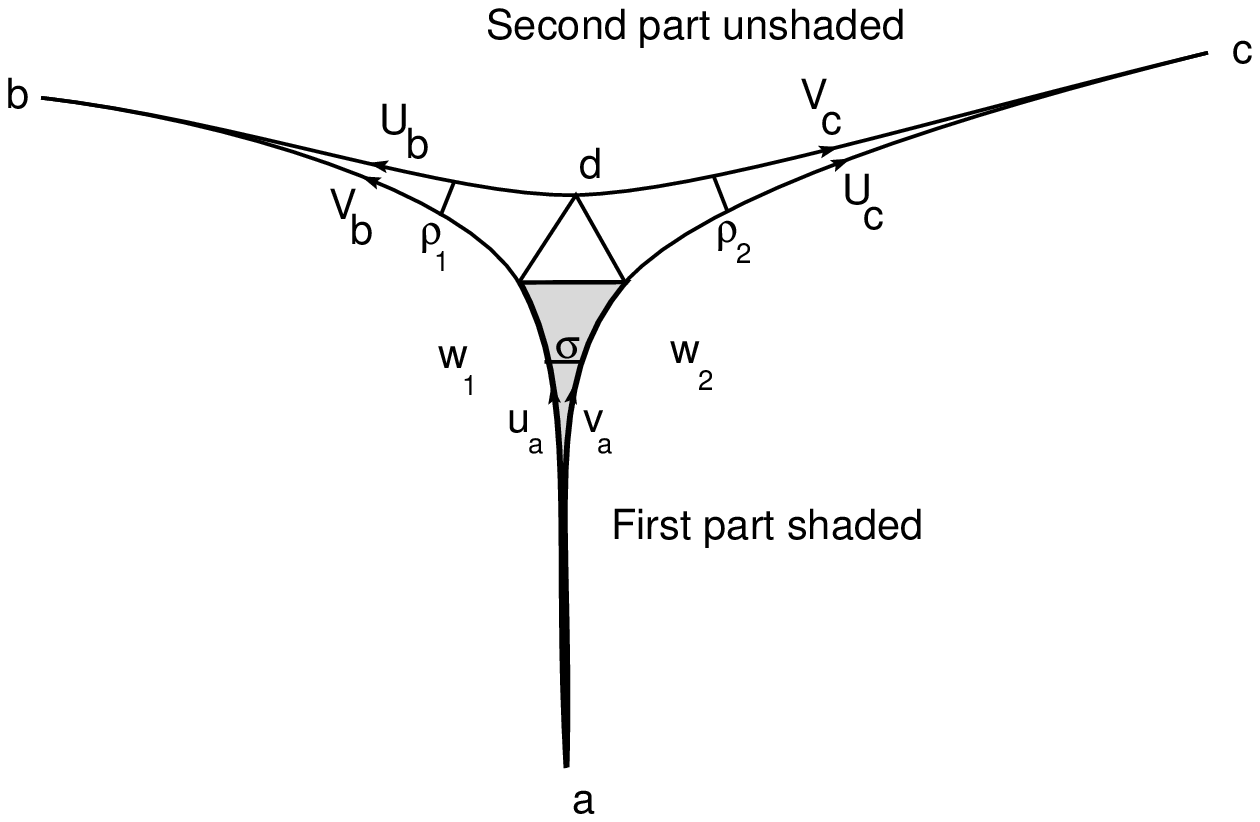}{This shows the path in the automaton, broken into a first
part which is in the automaton $FRD$ and a second part which is in two
copies of the automaton $FRD^R$.}

So the accepting path of $(w_1,w_2)$ through $NGP$ will be in two parts,
the first $(u_a,v_a)$ and the second $(V_b,U_c)$. A picture of this is
shown in Figure~\ref{complex}. Furthermore,
we have either $l(u_a) = l(v_a)$ or $l(u_a) = l(v_a)+1$, depending on
whether the perimeter of the geodesic triangle which has $w_1$ and
$w_2$ as two of its sides is even or odd.
In the first case, it is possible for $(V_b,U_c)$ to be empty,
which occurs when the vertices $b$ and $c$ of the geodesic triangle
coincide.

The first part of the
path through $NGP$ is simply a path through $FRD$,
ending at a state $\sigma \in \mathcal{A}(FRD)$.
In Figure~\ref{complex}, an intermediate state is also denoted by $\sigma$.

The second part corresponds to the two paths
$(U_b,V_b)$ and $(U_b,U_c)$ through $FRD^R$.
These paths must end at the unique accepting state of
$FRD^R$. This part of $NGP$ is non-deterministic,
because we need to quantify over their second variables as described in
Section~\ref{notation}. The initial states $\pi_1 = \{\sigma_1\}$ and
$\pi_2 = \{\sigma_2\}$ of these two paths through $FRD^R$ must be such that
$(\sigma,\sigma_1,\sigma_2) \in \mathcal{A}(FRD^3)$. This is equivalent to
$(u_a,v_a,U_b^R,V_b^R,U_c^R,V_c^R) \in L(FRD^3)$. In Figure~\ref{complex}
intermediate states are denoted by $(\rho_1,\rho_2)$.

In our implementation of $NGP$, we prefer to avoid
$\varepsilon$-transitions, and so the non-deterministic
jump from the first to the second part of the path is combined with the
first transitions in the second part of the path. In the case where
there is a padding symbol, the last transition in the first part of
the path is combined with the first transitions in the second part.
An advantage of this is that we can
eliminate the use of the padding symbol in the middle of a word, which
can otherwise be quite troublesome to deal with (in terms of writing
special code to take the unnecessary padding symbols into account).

The jump also introduces a large amount of non-determinism into $NGP$.

The states of $NGP$ are triples $(\sigma,\rho_1,\rho_2)$, where
$\sigma \in \mathcal{S}(FRD) \cup \{\infty\}$
and $\rho_1, \rho_2 \in \mathcal{S}(FRD^R) \cup \{0,\infty\}$.
For each such state either
$\rho_1=\rho_2=0$ and $\sigma \neq \infty$. or
$\sigma=\infty$ and $\rho_1 \neq 0 \neq \rho_2$.
Informally, $0$ and $\infty$ as just introduced have the following
significance.
In the course of accepting a string, the three components $\sigma$,
$\rho_1$ and $\rho_2$ each have to pass through $FRD$ exactly once.
More precisely,
the component $\sigma$ passes once through $FRD$ during the first part
of the path  and the components $\rho_1$
and $\rho_2$ pass through $FRD^R$ once during the second part of the
path, each time moving from an
initial state to an accept state of $FRD$ or $FRD^R$ respectively.
It is convenient to assign the name $\infty$ to the state 
$\sigma$ after it has completed its passage through $FRD$. We need to attach
names to the states $\rho_i$ ($i=1,2$) during the first part of
the path, and we attach the name $0$ to remind us that $\rho_i$
has not yet started its passage through $FRD^R$. When a padding symbol
is read in $w_i$, the state $\rho_i$ is set to $\infty$ to remind us
that the passage of $\rho_i$ through $FRD^R$ is now complete.
In other words, the state is set to $\infty$
the next move after arriving at the vertex of the triangle.
We do not allow $\rho_1=\rho_2=\infty$, because we will stop if we reach
both vertices $b$ and $c$ simultaneously.

We can save space by storing the triple $(\sigma,0,0)$ as a pair
$(\sigma,0)$ and the triple $(\infty,\rho_1,\rho_2)$ as
a pair $(\rho_1,\rho_2)$. It is easy to see that this captures all the
information.
In this discussion, we continue to use more revealing
triples, rather than more concise pairs.

The unique initial state is
$(\sigma_0(FRD),0,0)$. The accepting states are $(\infty,\rho_1,\rho_2)$ where
either $\rho_1=\infty$ or $\rho_1$ lies in $\mathcal{A}(FRD^R)$
and the same is true for $\rho_2$. The reader may like to be reminded
that a state of $\mathcal{A}(FRD^R)$ is an accept state of $FRD^R$, and that
this is a subset of the set of states of $FRD$ which contains the initial state
$\sigma_0(FRD)$.

There is also another kind of accept state, corresponding to the situation
$b=c$, that is, that there are two different geodesics from $a$ to $b=c$.
This means that $w_1=_G w_2$.
Since we are dealing with short-lex geodesics, $w_1$ will be short-lex
from $a$ to $b$ and $w_2$ will be the reverse of a short-lex geodesic
from $c$ to $a$. Such an accept state has the form $(\sigma,0)$ where
$(\sigma,\sigma_0(FRD),\sigma_0(FRD)) \in \mathcal{A}(FRD^3)$.

There are other degenerate situations, but the others are all covered by
the main description, as the reader can easily verify.

There are three types of transitions of $NGP$, which we shall now
describe.  In general, we denote the label of such a transition by
$(x_{ab},x_{ac})$, where $x_{ab},x_{ac} \in A^\dagger$. (The idea is
that $x_{ab}$ represents a variable generator in the path from $a$
to $b$.) The first two types of transition correspond to the
transitions of the first and second parts of the accepting path of
$(w_1,w_2)$ through $NGP$, and the third type of transition corresponds
to a jump from the first to the second part.

Transitions of the first type have the form $(\sigma,0,0) \rightarrow
(\tau,0,0)$ where there is a transition $\sigma \rightarrow \tau$ of
$FRD$ with the same label. However, we must have $x_{ab},x_{ac} \in A$;
that is, $x_{ac}$ is not allowed to be the padding symbol $\$$. (Any
such transition $(x,\$)$ of $FRD$
will be absorbed into the jump, and combined with the
first transitions of the two copies of $FRD^R$ after the jump,
which are bound to have labels
of the form $(y,\$)$ as we see from Figure~\ref{intv}. Recall that we
do not want the padding symbol to occur in the middle of either of the
words $w_1$, $w_2$.  The strategy explained here avoids that danger.)

Transitions of the second type have the form $(\infty,\pi_1,\pi_2) \rightarrow
(\infty,\rho_1,\rho_2)$.  They occur whenever there exist $x_{db},x_{dc}
\in A^\dagger$ and transitions $\pi_1 \rightarrow \rho_1$ and $\pi_2
\rightarrow \rho_2$ of $FRD^R$ with labels $(x_{db},x_{ab})$ and
$(x_{ac},x_{dc})$, respectively.  There is the further restriction that
$x_{ab}=\$$ if and only if $x_{db}=\$$.  This means that our path has
previously arrived at the vertex $b$---see Figure~\ref{complex}. In
this case we have $\rho_1 = \infty$ and either $\pi_1 = \infty$ or $\pi_1
\in \mathcal{A}(FRD^R)$. Similarly, $x_{ac}=\$$ if and only if
$x_{dc}=\$$. In this case $\rho_2 = \infty$ and either $\pi_2 = \infty$ or
$\pi_2 \in \mathcal{A}(FRD^R)$.  As explained above, we cannot have
$x_{ab}=x_{ac} = \$$.  In other words, padding symbols occur only at
the end of at most one of the words $w_1,w_2$.

The transitions of the third type are jumps from $(\sigma,0,0)$ to
$(\infty,\rho_1,\rho_2)$. These are of two subtypes, depending on whether the
geodesic triangle defined by the accepting paths that pass through them
has perimeter of even or odd length.

Those of even perimeter subtype occur when there exist $x_{db},x_{dc}
\in A$ and initial states $\pi_1 = \{\sigma_1\}$, $\pi_2 = \{\sigma_2\}$
of $FRD^R$ with the property that $(\sigma,\sigma_1,\sigma_2) \in
\mathcal{A}(FRD^3)$.  Furthermore there are transitions $\pi_1 \rightarrow
\rho_1$ and $\pi_2 \rightarrow \rho_2$ with labels $(x_{db},x_{ab})$
and $(x_{ac},x_{dc})$ respectively.  There is the further restriction
that $x_{ab}=\$$ if and only if $x_{db}=\$$, and in this case we have
$\rho_1 = \infty$ and $\pi_1 \in \mathcal{A}(FRD^R)$.  Similarly,
$x_{ac}=\$$ if and only if $x_{dc}=\$$, and in this case $\rho_2 =
\infty$ and $\pi_2 \in \mathcal{A}(FRD^R)$.

Those of the odd perimeter subtype arise only for $x_{ab},x_{ac} \in A$,
and they occur when there is a transition $\sigma \rightarrow \sigma'$
with label $(x_{ab},\$)$ of $FRD$. Furthermore, there exists $x_{dc}
\in A$ and initial states $\pi_1 = \{\sigma_1\}$, $\pi_2 = \{\sigma_2\}$
of $FRD^R$ with the property that $(\sigma',\sigma_1,\sigma_2) \in
\mathcal{A}(FRD^3)$.  Also, there are transitions $\pi_1 \rightarrow
\rho_1$ and $\pi_2 \rightarrow \rho_2$ with labels $(x_{db},\$)$ and
$(x_{ac},\$)$ respectively.

\section{Examples} \label{examples}
In this final section, we describe the performance of these algorithms
on the following four examples.
\begin{eqnarray*}
G_1 &=&\langle a,b,c,d \myvert a^{-1}b^{-1}abc^{-1}d^{-1}cd = 1 \rangle,\\
G_2 &=& \langle a,b \myvert a^2 = b^3 = (ab)^7 = 1 \rangle,\\
G_3 &=& \langle a,b \myvert (b^{-1}a^2ba^{-3})^2 = 1 \rangle \mbox{ and }
\end{eqnarray*}
\begin{eqnarray*}
\lefteqn{G_4 = \langle a,b,c,d,e,f \myvert a^4=b^4=c^4=d^4=e^4=f^4=}\\
&&aba^{-1}e = bcb^{-1}f = cdc^{-1}a = ded^{-1}b = efe^{-1}c = faf^{-1}d=1
\rangle.
\end{eqnarray*}

Of these, $G_1$ is a surface group of a genus two torus, $G_2$ is the
$(2,3,7)$-von Dyck group, $G_3$ is obtained from one of the
well-known family of non-Hopfian Baumslag-Solitar groups by squaring
the single relator, and $G_4$ is the symmetry group of a certain
tessellation by dodecahedra of hyperbolic 3-space as featured in the
video `Not Knot' (\cite{NotKnot}). They are all word-hyperbolic groups.

For all four the verification of hyperbolicity, as described in
Section~\ref{verifying},
was relatively easy, with the first three examples completing in a few seconds,
and $G_4$ requiring about 20 seconds cpu-time. We present details of these
calculations in Table 1. The automata $WD_n$, $GE_n$, $GW_n$ and $T_n$
are as described in Section~\ref{verifying},
and the constants $\gamma$ and $\gamma^\prime$
are as defined before Proposition~\ref{short-lex thin implies thin}.
The notation $m \rightarrow n$ in the table means that an automaton
had $m$ states when it was first constructed, and it was then minimised
to an equivalent automaton with $n$ states. The last example demonstrates
the phenomenon that the automata involved are smaller when the data is
correct.

\begin{table}[htbp]
\caption{Verifying Hyperbolicity}
\begin{center}
\begin{tabular}{|c|c|c|c|c|c|c|c|}
\hline
Grp & $n$ & $\mathcal{S}(WD_n)$ & $\mathcal{S}(GE_n)$
& $\mathcal{S}(GW_n)$ & $\mathcal{S}(T_n)$ & $\gamma$ &
$\gamma^\prime$\\
\hline
$G_1$ & 1 & 33 & $121 \rightarrow 49$ & $49 \rightarrow 49$ & 265 &4&4\\
\hline
$G_2$ & 1 & 30 & $627 \rightarrow 92$ & $80 \rightarrow 52$ & 936 &&\\
& 2 & 32 & $664 \rightarrow 94$ & $78 \rightarrow 54$ & 769 &7&7\\
\hline
$G_3$ & 1 & 55 & $689 \rightarrow 136$ & $152 \rightarrow 96$ & 1270 &6&6\\
\hline
$G_4$ & 1 & 75 & $896 \rightarrow 284$ & $454 \rightarrow 409$ & 10635 &&\\
& 2 & 97 & $1135 \rightarrow 309$ & $443 \rightarrow 378$ & 12407 &&\\
& 3 & 103 & $1211 \rightarrow 318$ & $424 \rightarrow 63$ & 1713 &4&4\\
\hline
\end{tabular}
\end{center}
\end{table}

In Table 2, we present details of the calculation of the thinness constant
for short-lex geodesic hyperbolic triangles in the first three of the examples.
The set $\mathcal{D}$ and the automata
$FRD$, $FRD^3$, $NGP$ and $GP$ are as defined in
Section~\ref{finding} (where $NGP$ is the non-deterministic version
of $GP$).
The separate lines of data for each group represent successive attempts at
the computation, with the last line representing the correct data.
After each attempt, the automaton with language
$L = L(W) \times L(W^R) \setminus L(GP)$ was constructed and, when it
was nonempty, words $ (w_1,w_2) \in L$
were found and used to find an improved set $\mathcal{D}$.
The language $L$ was found to be empty after the final computation for each
group, thereby proving correctness of the data.

The behaviour of $G_3$, which is the most difficult example for which we
have successfully completed the calculations, is probably the best indicator
of the way in which more difficult examples are likely to behave.
For example, the largest and most memory intensive part of the
computation is the determinisation of $NGP$, and many parts of the
calculations are significantly more expensive on the earlier passes,
when the data is incorrect, than in the final correct stage.
These computations were carried out using a maximum of 256 megabytes
of core memory and about the same amount of swap space.

We have not yet been able to complete the calculations for $G_4$; indeed
we have not progressed further than the first construction of $NGP$,
which has several million states.
We will need more memory (probably more than a gigabyte) if we are to proceed
to construct the determinised version $GP$.

\begin{table}[htbp]
\caption{Finding the Constant of Hyperbolicity}
\begin{center}
\begin{tabular}{|c|c|c|c|c|c|c|}
\hline
Grp & $ \mathcal{D} $ &
$\mathcal{S}(FRD)$ & $\mathcal{A}(FRD^3)$
& $\mathcal{S}(NGP)$ & $\mathcal{S}(GP)$ & $\delta$\\
\hline
$G_1$ & 25 & 137 & 65785 & 12249 & $8049 \rightarrow$ 2185 & \\
& 49 & 169 & 65857 & 12281 & $5457 \rightarrow 625$ & 4 \\
\hline
$G_2$ & 104 & 1174 & 73822 & 89802 & $35824 \rightarrow 4904$ & \\
& 111 & 1199 & 74047 & 90450 & $31374 \rightarrow 1508$  & 7 \\
\hline
$G_3$ & 71 & 755 & 795436 & 274186 & $1872679\rightarrow 531434 $ & \\
& 195 & 1430 & 801745 & 280240 & $1695944\rightarrow 443570 $ & \\
& 241 & 1741 & 806923 & 284328 & $1237158\rightarrow 85044 $ & \\
& 257 & 1845 & 807136 & 284478 & $676645 \rightarrow 3803 $ & 8 \\
\hline
\end{tabular}
\end{center}
\end{table}
\FloatBarrier


\begin{thebibliography}{99}
\bibitem[Alonso et al, 1991]{ALO}
J. Alonso, T. Brady, D. Cooper, V. Ferlini, M. Lustig, M. Mihalik,
M. Shapiro and H. Short, {\it Notes on word-hyperbolic groups}, in
E. Ghys, A. Haefliger and A. Verjovsky, eds.,
Proceedings of the Conference {\it Group Theory from a Geometric Viewpoint}
held in I.C.T.P., Trieste, March 1990, World Scientific, Singapore, 1991.

\bibitem[Epstein {\it et al.}, 1992]{ECHLPT}
D.B.A. Epstein, J.W. Cannon,
D.F. Holt, S.V.F. Levy, M.S. Paterson, W. Thurston (1992).
{\it Word Processing in Groups\/}, AKPeters, Natick, Mass.

\bibitem[Holt, 1995]{KBMAG}Derek F. Holt,
{\it {\sf KBMAG}---Knuth-Bendix in Monoids and Automatic Groups}, software
package (1995), available by anonymous {\tt ftp} from
{\tt ftp.maths.warwick.ac.uk} in directory {\tt people/dfh/kbmag2}.

\bibitem[Holt, 1996]{HOLT}D.F. Holt,
{\it The Warwick automatic groups software},
in {\it Geometrical and Computational Perspectives on Infinite Groups}
DIMACS Series in Discrete Mathematics and Theoretical Computer Science,
vol. 25, ed. G. Baumslag et. al., 1996, 69--82.

\bibitem[Not Knot]{NotKnot}Charlie Gunn and Delle Maxwell,
{\it Not Knot}, video ISBN 0-86720-240-8, published by AKPeters,
Natick, Mass., with a booklet by D.B.A. Epstein and C. Gunn.

\bibitem [Papasoglu, 1994]{PAP}, P. Papasoglu,
{\it Strongly geodesically automatic groups
are hyperbolic}, (Warwick preprint) (1994) (to appear).

\bibitem[Wake\-field, 1997]{WAK} P. Wakefield,
{\it Procedures for Automatic groups},
Ph.D. Thesis, University of Newcastle upon Tyne, 1997.
\end{thebibliography}
\end{document}